\documentclass[12pt,leqno]{article}
\usepackage{amssymb}
\usepackage{amscd}
\usepackage{amsmath}
\usepackage{amsthm}
\usepackage{graphics}
\usepackage{mathrsfs}
\def\overset#1#2{{\mathrel{\mathop {{#2}_{}}\limits^{#1}}}}
\def\underset#1#2{{\mathrel{\mathop {{}_{} {#2}}\limits_{{#1}_{}}}}}
\def\upplim_#1{\underset{#1}{\overline\lim}\;}
\def\lowlim_#1{\underset{#1}{\underline\lim}\;}
\newtheorem{cond}[equation]{\indent \rm Condition}
\newtheorem{cor}[equation]{Corollary}
\newtheorem{defn}[equation]{\indent{\it Definition}\rm }

\newtheorem{lem}[equation]{Lemma}

\newtheorem{rmk}[equation]{\indent \rm {\it Remark}}

\newtheorem{thm}[equation]{Theorem}

\newcommand{\C}{{\mathbf{C}}}

\newcommand{\del}{{\partial}}

\newcommand{\iso}{\cong}

\newcommand{\sI}{\mathscr{I}}

\newcommand{\lto}{\longrightarrow}

\newcommand{\N}{\mathbf{N}}

\renewcommand{\O}{{\mathcal{O}}}

\newcommand{\PD}{\mathrm{P}\Delta}
\newcommand{\pnc}{{\mathbf{P}^n(\mathbf{C})}}

\newcommand{\Z}{\mathbf{Z}}
\newenvironment{pf}{\par{\it Proof. }}{\qed}
\numberwithin{equation}{section}

\title{
A Scalar Associated with
 the Inverse of
 Some
 Abelian \\ Integrals 
 and  a Ramified Riemann Domain 
}
\author{
Junjiro Noguchi\thanks{
  Research supported in part by Grant-in-Aid
 for Scientific Research (B) 23340029. \hfill\break
Keywords: Abelian integral; open Riemann surface; Levi problem;
 Riemann domain; Stein manifold.}
}
\begin{document}
\setlength{\baselineskip}{16pt}
\maketitle
\begin{abstract}
We introduce a positive scalar function $\rho(a, \Omega)$ for a
domain $\Omega$ of a complex manifold $X$ with a global
 holomorphic frame of the cotangent bundle by closed
Abelian differentials,
which 
is an analogue of Hartogs' radius.
We prove an {\em estimate
of Cartan--Thullen type with $\rho(a, \Omega)$}
for holomorphically convex hulls of compact subsets.
In one dimensional case,  we apply the obtained estimate
 of $\rho(a, \Omega)$ 
to give a new proof of Behnke-Stein's Theorem for the Steiness
of open Riemann surfaces.
We then extend the idea to deal with the problem to generalize
Oka's Theorem (IX) for ramified Riemann domains over $\C^n$.
We obtain some geometric conditions in terms of $\rho(a, X)$
which imply the validity of the Levi problem (Hartogs' inverse problem)
for a finitely sheeted Riemann domain over $\C^n$.
\end{abstract}

\section{Introduction and main results}

\subsection{Introduction}
In 1943 K. Oka wrote a manuscript in Japanese, solving
affirmatively the Levi problem (Hartogs' inverse problem)
 for unramified Riemann domains
over complex number space $\C^n$ of arbitrary dimension
$n \geq 2$,\footnote{
This fact was written twice in the introductions of his two papers,
\cite{ok8} and \cite{ok9}: The manuscript was written as a
research report {\em dated 12 Dec.\ 1943},
 sent to Teiji Takagi, then Professor at
the Imperial University of Tokyo, and now one can find it
in \cite{oka-bunko}.
} and 
in 1953 he published Oka IX \cite{ok9} to solve it by making use of
his First Coherence Theorem proved in Oka VII \cite{ok7}\,\footnote{
It is noted that Oka VII \cite{ok7} is different to
his original, Oka VII in \cite{ok}; therefore, there are two
versions of Oka VII. The English translation of Oka VII in \cite{oksp}
was taken from the latter, but unfortunately in \cite{oksp}
all the records of the received dates of the papers were deleted.
};
there, he put a special emphasis on the difficulties of the ramified case
(see \cite{ok9}, Introduction 2 and \S23).
H.~Grauert also emphasized the problem to generalize Oka's Theorem (IX)
 to the casse of  ramified Riemann domains in his lecture at
  OKA 100 Conference Kyoto/Nara 2001.
Oka's Theorem (IX) was generalized for unramified Riemann domains
over complex projective $n$-space $\pnc$ by R. Fujita \cite{fuj}
and A. Takeuchi \cite{tak}.
Later, H. Grauert \cite{nar} gave a counter-example to the problem
for ramified Riemann domains over $\pnc$, and 
J.E. Forn{\ae}ss \cite{forn} gave a counter-example
 to it over $\C^n$.
Therefore, it is natural to look for  geometric conditions
which imply the validity of the Levi problem for ramified
Riemann domains.

Under a geometric condition (Cond\,A, \ref{conda}) on a
complex manifold $X$,  we introduce 
a new scalar function
$\rho(a, \Omega) (>0)$ for a subdomain $\Omega \subset X$,
which is an analogue of the boundary distance
function in the unramified case (cf.\ Remark \ref{hradius} (i)).
We prove an {\em estimate of Cartan-Thullen type} (\cite{ca-th})
for the holomorphically convex hull $\hat{K}_\Omega$
of a compact subset $K \Subset \Omega$ with $\rho(a, \Omega)$
(see Theorem \ref{ct}).

In one dimensional case, by making use of $\rho(a, \Omega)$
we give a {\em new proof of Behnke-Stein's Theorem}: Every open Riemann
surface is Stein.
In the known methods one uses a generalization of the
Cauchy kernel or some functional analytic method (cf.\
Behnke-Stein \cite{be-st}, Kusunoki \cite{kus}, Forster \cite{fors}, etc.).
Here we use Oka's J\^oku-Ik\^o combined with
Grauert's finiteness theorem, which is now a rather easy
result by a simplification of the proof,
 particularly in $1$-dimensional case (see \S\ref{rs}).
We see here how the scalar $\rho(a, \Omega)$ works well in this case.

Now, let $\pi:X \to \C^n$ be a Riemann domain, possibly ramified,
such that $X$ satisfies Cond\,A. Then, we prove that
a domain $\Omega \Subset X$ is a {\em domain of holomorphy\footnote{
The notion of ``domain of holomorphy'' for $\Omega \Subset X$ is
defined as usual (cf., e.g., \cite{ho}).
} if and only
if $\Omega$ is holomorphically convex} (see Theorem \ref{relcomp}).
Moreover, if $X$ is exhausted by a continuous family of
relatively compact domains of holomorphy,
then $X$ is Stein (see Theorem \ref{exh1}).

We next consider a boundary condition (Cond\,B, \ref{condb})
with $\rho(a, X)$.
We assume that $X$ satisfies Cond\,A and that
 $X ~\overset{\pi}{\to}~ \C^n$ satisfies
Cond\,B and is finitely sheeted. We prove that
{\em if $X$ is locally Stein over $\C^n$, then $X$ is Stein}
(see Theorem \ref{levi2}).

We give the proofs in \S2. In \S3 we will discuss
some examples and properties of $\rho(a, X)$.

\smallskip
{\it Acknowledgment.}  The author is very grateful to
Professor Makoto Abe for interesting discussions
on the present theme.

\subsection{Main results}

\subsubsection{Scalar $\rho(a, \Omega)$}\label{scr}
Let $X$ be a connected complex manifold of dimension $n$
 with holomorphic cotangent
bundle $\mathbf{T}(X)^*$. We assume:

\begin{cond}[Cond\,A]\rm
\label{conda}
There exists a global frame $\omega=(\omega^1, \ldots, \omega^n)$
of $\mathbf{T}(X)^*$ over $X$ such that $d\omega^j=0$,
$1 \leq j \leq n$.
\end{cond}

Let $\Omega \subset X$ be a subdomain.
With Cond\,A we consider
 an Abelian integral (a path integral)
of $\omega$ in $\Omega$ from $a \in \Omega$:
\begin{equation}
\label{abint}
\alpha: x \in \Omega \lto  \zeta= (\zeta^j)=
\left(\int_a^x \omega^1, \ldots, \int_a^x \omega^n
\right)  \in \C^n.
\end{equation}
We denote by $\PD=\prod_{j=1}^n \{|\zeta^j|< 1\}$ the 
unit polydisk  of
$\C^n$ with center at $0$ and and set
$$
\rho\PD=\prod_{j=1}^n \{|\zeta^j|< \rho\}
$$
for $\rho>0$.
Then, $\alpha(x)=\zeta$ has the inverse $\phi_{a,\rho_0}(\zeta)=x$ on
a small polydisk $\rho_0\PD$:
\begin{equation}
\label{locbih}
 \phi_{a,\rho_0}: \rho_0 \PD \lto U_0=\phi_{a,\rho_0}
(\rho_0 \PD) \subset \Omega.
\end{equation}
Then we extend analytically $\phi_{a,\rho_0}$
to $\phi_{a,\rho}: \rho \PD \to X$, $\rho \geq \rho_0$,
 as much as possible, and set
\begin{equation}
\label{convrad}
\rho(a, \Omega) =\sup\{\rho>0 : \exists \phi_{a,\rho}: \rho \PD
\to X, ~ \phi_{a,\rho}(\rho \PD) \subset  \Omega\}.
\end{equation}
Then we have the inverse of the Abelian integral $\alpha$
on the polydisk of the maximal radius
\begin{equation}
\label{invabint}
\phi_a : \rho(a, \Omega) \PD \lto \Omega.
\end{equation}
To be precise, we should write
\begin{equation}
\label{rho}
 \rho(a, \Omega)=\rho(a, \omega, \Omega)=\rho(a, \PD, \omega, \Omega),
\end{equation}
but unless confusion occurs,
 we use $\rho(a, \Omega)$ for notational simplicity.

We immediately see that (cf.\ \S\ref{sclr})
\begin{enumerate}
\setlength{\itemsep}{-2pt}
\item
$\rho(a, \Omega)$ is continuous;
\item
$\rho(a, \Omega) \leq \inf\{|v|_\omega\,{:}\, v \in
\mathbf{T}(X)_a ,\, F_\Omega(v)=1\}$,
where $F_\Omega$ denotes the Kobayashi hyperbolic infinitesimal form
of $\Omega$, and $|v|_\omega=\max_j |\omega^j(v)|$,
 the maximum norm of $v$ with respect to $\omega=(\omega^j)$.
\end{enumerate}

For a subset $A \subset \Omega$ we write
\[
 \rho(A, \Omega)=\inf\{ \rho(a, \Omega): a \in A\}.
\]

For a compact subset $K \Subset \Omega$ we denote by $\hat{K}_\Omega$
the holomorphically convex hull of $K$ defined by
\[
 \hat{K}_\Omega=\left\{x \in \Omega: |f(x)| \leq \max_K
|f|, ~ \forall f \in \O(\Omega)\right\},
\]
where $\O(\Omega)$ is the set of all holomorphic functions
on $\Omega$.
If $\hat{K}_\Omega \Subset \Omega$ for every $K \Subset \Omega$,
$\Omega$ is called a holomorphically convex domain.

\begin{defn}\rm
For a relatively compact subdomain $\Omega \Subset X$ of a complex
manifold $X$
we may naturally define the notion of {\em domain of holomorphy}: i.e.,
there is no point $b \in \del \Omega$ such that there are
a connected neighborhood $U$ of $b$ in $X$
and a non-empty open subset $V \subset U \cap \Omega$
satisfying that for every $f \in \O(\Omega)$
there exists $g \in \O(U)$ with $f|_V=g|_V$.
\end{defn}

The following theorem of Cartan-Thullen type (cf.\ \cite{ca-th})
is our first main result.
\begin{thm}
\label{ct}
Let $X$ be a complex manifold satisfying Cond\,A.
Let $\Omega \Subset X$ be a relatively compact domain of holomorphy,
let $K \Subset \Omega$ be a compact subset,
and let $f\in \O(\Omega)$.
Assume that
\[
|f(a)| \leq \rho(a, \Omega), \quad \forall a \in K.
\]
Then we have
\begin{equation}
\label{ct1}
|f(a)| \leq \rho(a, \Omega), \quad \forall a \in \hat{K}_\Omega.
\end{equation}
In particular, we have
\begin{equation}
\label{ct2}
\rho(K, \Omega)=\rho(\hat{K}_\Omega, \Omega).
\end{equation}
\end{thm}

\begin{cor}
\label{dhol}
Let $\Omega \Subset X$ be a domain of a complex manifold $X$,
satisfying Cond A.
Then, $\Omega$ is a domain of holomorphy if and only if $\Omega$
is holomorphically convex.
\end{cor}

\subsubsection{Behnke-Stein's Theorem for open Riemann surfaces}\label{rs}

We apply the scalar $\rho(a, \Omega)$ introduced above
to give a new proof of Behnke-Stein's Theorem for the Steiness
of open Riemann surfaces, which is one of the most basic facts
in the theory of Riemann surfaces: Here, we do not use
the Cauchy kernel generalized on a Riemann surface
(cf.\ \cite{be-st}, \cite{kus}), nor a functional analytic method
(cf., e.g., \cite{fors}),
 but use Oka's J\^oku-Ik\^o together with
Grauert's finiteness theorem, which is now a rather easy
result, particularly in $1$-dimensional case.
This is the very difference of our new proof to the known ones.

To be precise, we recall the definition of Stein manifold:
\begin{defn}\label{stdef}
\rm
A  complex manifold $M$ of pure
dimension $n$ is called a Stein manifold if the following
Stein conditions are satisfied:
\begin{enumerate}
\setlength{\itemsep}{-2pt}
\item
$M$ satisfies the second countability axiom.
\item
For distinct points $p,q \in M$ there is an $f \in \O(M)$
with $f(p)\not=f(q)$.
\item
For every $p \in M$ there are $f_j \in \O(M)$, $1 \leq j \leq n$,
such that $df_1(p) \wedge \cdots \wedge df_n (p) \not =0$.
\item
$M$ is holomorphically convex.
\end{enumerate}
\end{defn}

We will rely on
the following H. Grauert's Finiteness Theorem in $1$-dimensional case,
 which is now  a rather easy consequence of Oka--Cartan's Fundamental
Theorem, particularly in $1$-dimensional case,
 thanks to a very simplified proof of L. Schwartz's
Finiteness Theorem based on the idea of
 Demailly's Lecture Notes \cite{dem},
Chap.\ IX (cf.\ \cite{nog1}, \S7.3 for the present form):

\smallskip
{\bf L. Schwartz' Finiteness Theorem}. 
{\em Let $E$ be a Fr\'echet space and let $F$ be a Baire vector space.
Let $A: E \to F$ be a continuous linear surjection, 
and let $B: E \to F$ be a completely continuous linear map.
Then, $(A+B)(E)$ is closed and the cokernel 
$\mathrm{Coker}(A+B)$ is finite dimensional.}

Here, a Baire space is a topological space such that Baire's category theorem holds.
The statement above is slightly generalized than the original one,
in which $F$ is also assumed to be Fr\'echet (cf.\ L. Schwartz \cite{sch},
Serre \cite{ser}, Bers \cite{ber},
 Grauert-Remmert \cite{gr2}, Demailly \cite{dem}).

\smallskip
{\bf Grauert's Theorem in dimension 1}.
{\em Let $X$ be a Riemann surface, and let 
$\Omega \Subset X$ be a relatively compact subdomain. Then,}
\begin{equation}
\label{gra}
 \dim H^1(\Omega, \O_\Omega) < \infty.
\end{equation}

Here, $\O_\Omega$ denotes the sheaf of germs of holomorphic
functions over $\Omega$. In case $\Omega (=X)$ itself is compact,
this theorem reduces to Cartan--Serre's in dimension $1$.

{\bf N.B.} It is the very idea of Grauert to claim only the finite
dimensionality, weaker than a posteriori statement,
$H^1(\Omega, \O_\Omega)=0$: It makes the proof considerably easy.

\smallskip
By making use of this theorem we prove an intermediate result:
\begin{lem}
\label{relcpt}
Every relatively compact domain $\Omega$ of $X$ is Stein.
\end{lem}

Let $\Omega \Subset \tilde{\Omega} \Subset X$ be subdomains
of an open Riemann surface $X$.
Since $\tilde{\Omega}$ is Stein by Lemma \ref{relcpt}
 and $H^2(\tilde{\Omega}, \Z)=0$,
we see by the Oka Principle that the line bundle
of holomorphic 1-forms over $\tilde{\Omega}$
is trivial, and so we have:
\begin{cor}
There exists a holomorphic
$1$-form $\omega$ on $\tilde{\Omega}$ without zeros.
\end{cor}

By making use of $\omega$ above we define $\rho(a, \Omega)$
as in \eqref{convrad} with $X=\tilde{\Omega}$.

Applying Oka's J\^oku-Ik\^o combined  with $\rho(a, \Omega)$,
 we give the proofs of
the following approximations of Runge type:

\begin{lem}
\label{rung1}
Let $\Omega'$ be a domain such that $\Omega \Subset \Omega' \Subset
\tilde{\Omega}$. Assume that
\begin{equation}
\label{bcond}
\max_{b \in \del\Omega} \rho(b, \Omega') < \rho(K, \Omega).
\end{equation}
Then, every $f \in \O(\Omega)$ can be approximated uniformly
on $K$ by elements of $\O(\Omega')$.
\end{lem}

\begin{thm}
\label{rung3}
Assume that no component of $\tilde{\Omega} \setminus \bar{\Omega}$
is relatively compact in $\tilde{\Omega}$.
Then, every $f \in \O(\Omega)$ can be approximated
uniformly on compact subsets of $\Omega$ by elements
of $\O(\tilde{\Omega})$.
\end{thm}

Finally we give a proof of
\begin{thm}[Behnke-Stein \cite{be-st}]
\label{bs}
Every open Riemann surface $X$ is Stein.
\end{thm}

\subsubsection{Riemann domains}\label{subsrdom}
Let $X$ be a complex manifold, and let $\pi:X \to \C^n$
(resp.\ $\pnc$) be a holomorphic map.
\begin{defn}
\label{rdom}\rm
We call $\pi: X \to \C^n$ (resp.\ $\pnc$) a {\em Riemann domain}
 (over $\C^n$ (resp.\  $\pnc$))
 if every fiber  $\pi^{-1}z$ with $z \in \C^n$ (resp.\ $\pnc$) is discrete;
if $d\pi$ has the maximal rank everywhere,
it is called an {\em unramified}
Riemann domain (over $\C^n$ (resp.\ $\pnc$)).
A Riemann domain which is not unramified, is called
 a {\em ramified Riemann domain}.
If the cardinality of $\pi^{-1}z$ is bounded in
 $z \in \C^n$ (resp.\ $\pnc$),
then we say that $\pi: X \to \C^n$ (resp.\ $\pnc$)
 is {\em finitely sheeted} or {\em $k$-sheeted} with
the maximum $k$ of the cardinalities  of $\pi^{-1}z$ ($z \in \C^n$
(resp.\ $\pnc$)).
\end{defn}

If $\pi: X \to \C^n$ (resp.\ $\pnc$) is a Riemann domain,
then the pull-back of Euclidean metric (resp.\ Fubini-Study metric)
by $\pi$ is a degenerate (pseudo-)hermitian metric on $X$,
which leads a distance function on $X$; hence, $X$
satisfies the second countability axiom.

Note that unramified Riemann domains over $\C^n$ naturally
 satisfy Cond\,A.

We have:
\begin{thm}
\label{relcomp}
Let $\pi:X \to \C^n$ be a Riemann domain such that
$X$ satisfies Cond\,A. 
\begin{enumerate}
\setlength{\itemsep}{-2pt}
\item
Let $\Omega \Subset X$ be a subdomain.
Then, $\Omega$ is a domain of holomorphy if and only if
$\Omega$ is Stein.
\item
If $X$ is Stein,
then $-\log \rho(a, X)$ is either identically $-\infty$,
or continuous plurisubharmonic.
\end{enumerate}
\end{thm}

\begin{defn}[Locally Stein]\rm
\label{lstein}
\begin{enumerate}
\item
Let $X$ be a complex manifold.
We say that a subdomain $\Omega \Subset X$ is {\em locally Stein} if
for every $a \in \bar{\Omega}$ (the topological closure)
there is a neighborhood $U$ of $a$ in $X$ such that
$\Omega \cap U$ is Stein.
\item
Let $\pi:X \to \C^n$ be a Riemann domain.
If for every point $z \in \C^n$ there is a
neighborhood $V$ of $z$
such that $\pi^{-1}V$ is Stein, $X$ is said to be
 {\em locally Stein over $\C^n$} (cf.\ \cite{forn}).
\end{enumerate}
\end{defn}
In general, the Levi problem is the one to asks if a locally Stein domain (over $\C^n$)
 is Stein.

\begin{rmk}\rm
\label{levi1}
The following statement is a direct consequence of 
Elencwajg \cite{el}, Th\'eor\`eme II combined with
Andreotti-Narasimhan \cite{an}, Lemma 5:

\begin{thm}\label{el-annar}
Let $\pi:X \to \C^n$ be a Riemann domain,
and let $\Omega \Subset X$ be a subdomain.
If $\Omega$ is locally Stein, then $\Omega$ is a
Stein manifold.
\end{thm}

Therefore the Levi problem for a ramified Riemann domain
$X\,\overset{\pi}{\to}\,\C^n$ is essentially at the ``{\em infinity}''
 of $X$.
\end{rmk}

\begin{defn}\label{cex}\rm
Let $X$ be a complex manifold in general.
 A family $\{\Omega_t\}_{0 \leq t \leq 1}$
of subdomains $\Omega_t$ of $X$ is called a
 {\em continuous exhaustion family of subdomains of $X$} if the
following conditions are satisfied:
\begin{enumerate}
\item
$\Omega_t \Subset \Omega_s \Subset \Omega_1=X$ for $0 \leq t < s <1$,
\item
$\bigcup_{t <s} \Omega_t = \Omega_s $ for $0  < s \leq 1$,
\item
$\del \Omega_t =\bigcap_{s>t} \overline{{\Omega}_s \setminus 
\overline{\Omega}_t}$
for $0 \leq t <1$.
\end{enumerate}
\end{defn}

\begin{thm}\label{exh1}
Let $\pi: X \to \C^n$ be a Riemann domain.
Assume that there is a continuous exhaustion family
$\{\Omega_t\}_{0 \leq t \leq 1}$ of subdomains of $X$ such that
for $0 \leq t <1$,
\begin{enumerate}
\item
 $\Omega_t$  satisfies Cond\,A,
\item
 $\Omega_t$ is a domain of holomorphy
 (or equivalently, Stein).
\end{enumerate}
Then, $X$ is Stein, and for any fixed $0 \leq t <1$
a holomorphic function $f \in \O(\Omega_t)$ can be
approximated uniformly on compact subsets by elements of
$\O(X)$.
\end{thm}

Let $\pi: X \to \C^n$ be a Riemann domain
such that $X$ satisfies Cond\,A
and let $\del X$ denote the ideal boundary of $X$ over $\C^n$
(called the accessible boundary in Fritzsche-Grauert
 \cite{f-g}, Chap.\ II \S9).
To deal with the total space $X$ we consider the following condition
which is a sort of {\em localization principle}:
\begin{cond}[Cond\,B]\rm
\label{condb}
\begin{enumerate}
\item
$\lim_{a \to \del X} \rho(a, X)=0$,
\item 
For every ideal boundary point $b \in \del X$
there are neighborhoods $V \Subset W$  of 
 $\pi(b)$  in $\C^n$
such that 
for the connected components $\widetilde{V}$ of 
$\pi^{-1}V$ and  $\widetilde{W}$ of $\pi^{-1}W$ with
 $\widetilde{V} \subset \widetilde{W}$,
which are elements of the defining filter of $b$,
\begin{equation}
\label{lp}
\rho(a, X) = \rho(a, \widetilde{W}), \quad \forall a \in
\widetilde{V}.
\end{equation}
\end{enumerate}
\end{cond}

For the Levi problem in case (ii) we prove:
\begin{thm}
\label{levi2}
Let $\pi: X \to \C^n$ be a finitely sheeted
Riemann domain. Assume that Cond\,A and Cond\,B
are satisfied.
If $X$ is locally Stein over $\C^n$,  $X$ is
a Stein manifold.
\end{thm}

\begin{rmk}\rm
Forn{\ae}ss' counter-example (\cite{forn})
 for the Levi problem in the ramified case
is a $2$-sheeted Riemann domain over $\C^n$.
\end{rmk}

\section{Proofs}
\subsection{Scalar $\rho(a, \Omega)$}\label{sclr}
Let $X$ be a complex manifold satisfying Cond\,A.
We here deal with some elementary properties of $\rho(a, \Omega)$
defined by \eqref{convrad} for a subdomain $\Omega \Subset X$.
We use the same notion as in \S\ref{scr}.
We identify $\rho \PD_0$ and $U_0$ in \eqref{locbih}.
For $b, c \in \rho_0 \PD$ we have
\[
 \rho(b, \Omega) \geq \rho(c, \Omega) - |b-c|,
\]
where $|b-c|$ denotes the maximum norm with respect to
the coordinate system $(\zeta^j) \in \rho_0 \PD$.
Thus,
\[
 \rho(c, \Omega) - \rho(b, \Omega) \leq |b-c|.
\]
Changing $b$ and $c$, we have the converse inequality, so that
\begin{equation}
\label{rcont}
|\rho(b, \Omega) - \rho(c, \Omega)| \leq |b-c|,
\quad b, c \in \rho_0 \PD \iso U_0.
\end{equation}
Therefore, $\rho(a, \Omega)$ is a continuous function in $a \in \Omega$.

Let
$v=\sum_{j=1}^n v^j \left(\frac{\del~}{\del \zeta^j}\right)_a \in
\mathbf{T}(\Omega)_a$
 be a holomorphic tangent at $a \in \Omega$.
Then,
$$
|v|_\omega=\max_{1 \leq j \leq n} |v^j|.
$$
With $|v|_\omega=1$ we have by the definition of the Kobayashi
hyperbolic infinitesimal metric $F_\Omega$
(cf.\ \cite{ko}, \cite{nw})
\[
 F_\Omega(v) \leq \frac{1}{\rho(a, \Omega)}.
\]
Therefore we have
\begin{equation}
\label{kh}
\rho(a, \Omega) \leq \inf_{v: F_\Omega(v)=1} |v|_\omega.
\end{equation}

Provided that $\del \Omega \not= \emptyset$, it immediately follows that
\begin{equation}
\label{b0}
\lim_{a \to \del\Omega} \rho(a, \Omega)=0.
\end{equation}

\begin{rmk}\rm
\label{hradius}
(i) We consider an {\em unramified} Riemann domain $\pi:X \to \C^n$.
Let $(z^1, \ldots, z^n)$ be the natural coordinate system of $\C^n$
and put $\omega=(\pi^* dz^j)$. Then the boundary distance function
$\delta_{\PD}(a, \del X)$ to the ideal boundary $\del X$
with respect to the unit polydisk $\PD$ is
defined as the supremum of such $r>0$ that $X$
is univalent onto $\pi(a)+r\PD$ in a neighborhood of $a$
(cf., e.g.,  \cite{ho}, \cite{nog1}).
Therefore, in this case we have that
\begin{equation}
\label{rd}
 \rho(a, X)=\delta_{\PD}(a, \del X).
\end{equation}
As for the difficulty to deal with the Levi problem for ramified Riemann
domains, K. Oka wrote in IX \cite{ok9}, \S23:
\begin{quote}
`` Pour le deuxi\`eme cas
les rayons de \textit{Hartogs} cessent de jouir
du r\^ole; ceci pr\'esente une difficult\'e qui m'apparait vraiment
 grande.'' 
\end{quote}
The above ``le deuxi\`eme cas'' is the ramified case.

(ii) For $X$ satisfying Cond\,A one can define {\em Hartogs' radius}
$\rho_n(a, X)$ as follows. Consider $\phi_{a, (r_j)}: \PD(r_j) \to X$ for
a polydisk $\PD(r_j)$ about $0$ with
a poly-radius $(r_1, \ldots, r_n)$ ($r_j>0$), which is
an inverse of $\alpha$ given by \eqref{abint}. Then, one defines
$\rho_n(a, X)$ as the supremum of such $r_n>0$; for other $j$,
it is similarly defined. Hartogs' radius $\rho_n(a, \Omega)$
is not necessarily continuous, but lower semi-continuous.
In the present paper, the scalar $\rho(a, X)$ defined under Cond\,A
plays the role of ``Hartogs' radius''.
\end{rmk}

\begin{rmk}\rm
Let $X$ be a complex manifold satisfying Cond\,A.
We see that {\em if $\rho(a_0, X)=\infty$ at a point $a_0 \in X$,
then $\phi_{a_0}: \C^n \to X$ is surjective, and
$\rho(a, X)\equiv \infty$ for $a \in X$.}
In fact, suppose that $\rho(a_0, X)=\infty$. Then, for
any $a \in X$ we take a path $C_{a}$ from $a_0$ to $a$
and set $\zeta=\alpha(a)$. By the definition,
$\phi_{a_0}(\zeta)=a$, and it follows that
$\rho(a, X)=\infty$.
Even if $\rho(a,\omega, X)=\infty$ (cf.\ \eqref{rho}),
``$\rho(a,\omega', X)< \infty$'' may happen for another choice of $\omega'$
(cf.\ \S\ref{smr}).
\end{rmk}

\subsection{Proof of Theorem \ref{ct}}
For $a \in \Omega$ we let
$$
\phi_a: \rho(a, \Omega) \PD \lto \Omega
$$
be as in \eqref{invabint}.
We take an arbitrary element $u \in \O(\Omega)$.
With a fixed positive number $s <1$ we set
\[
 L=\bigcup_{a \in K} \phi_a \left(s |f(a)|\, \overline{\PD}\right).
\]
Then it follows from the assumption that $L$ is a compact subset
of $\Omega$. Therefore there is an $M>0$ such that
$$
|u| <M  \hbox{ on } L.
$$
Let $\del_j$  be the dual vector fields
of $\omega^j$, $1 \leq j \leq n$, on $X$.
For a multi-index $\nu=(\nu_1, \ldots, \nu_n)$
with non-negative integers $\nu_j \in \Z^+$ we put
\begin{align*}
 \del^\nu &=\del_1^{\nu_1}
\cdots \del_n^{\nu_n} \\
|\nu| &= \nu_1+ \cdots +\nu_n,\\
\nu! &= \nu_1! \cdot \cdots \cdot \nu_n ! \, .
\end{align*}

By Cauchy's inequalities for $u \circ \phi_a$
on $s |f(a)|\, \overline{\PD}$ with $a \in K$ we have
\[
\frac{1}{\nu!} |\del^\nu u(a)| \cdot
 |s f(a)|^{|\nu|} \leq M, \quad \forall a \in K,
~ \forall \nu \in (\Z^+)^n.
\]
Note that $(\del^\nu u) \cdot f^{|\nu|} \in \O(\Omega)$.
By the definition of $\hat{K}_\Omega$,
\begin{equation}
\label{2.7}
\frac{1}{\nu!} |\del^\nu u(a)| \cdot
 |s f(a)|^{|\nu|} \leq M, \quad
\forall a \in \hat{K}_\Omega, ~ \forall \nu \in (\Z^+)^n.
\end{equation}

For $a \in \hat{K}_\Omega$
we consider the Taylor expansion of $u \circ \phi_a (\zeta)$ at $a$:
\begin{align}
\label{tay}
u \circ \phi_a (\zeta) &=\sum_{\nu \in (\Z^+)^n}
\frac{1}{\nu!} \del^\nu u (a) \zeta^\nu.
\end{align}
We infer from \eqref{2.7} that \eqref{tay} converges at least
on $s|f(a)|\,\PD$. Since $\Omega$ is a domain of holomorphy,
we have that $\rho(a, \Omega) \geq s |f(a)|$.
Letting $s \nearrow 1$, we deduce \eqref{ct1}.

By definition, $\rho(K, \Omega) \geq \rho(\hat{K}_\Omega, \Omega)$.
The converse is deduced by applying the result obtained above
for a constant function $f \equiv \rho(K, \Omega)$;
thus \eqref{ct2} follows. \qed

\smallskip
{\bf Proof of Corollary \ref{dhol}:}
 Assume that $\Omega \Subset X$ is a domain of holomorphy.
Let $K \Subset \Omega$. It follows from \eqref{ct2} that
$\hat{K}_\Omega \Subset \Omega$, and hence $\Omega$
is holomorphically convex.
The converse is clear. \qed

\begin{rmk}\rm
(i) Replacing $\PD$ by the unit ball $B$ with center at $0$,
 one may define similarly $\rho(a, \Omega)$.
Then Theorem \ref{ct} remains to hold.
Note that the union of all unitary rotations of
$\frac{1}{\sqrt{n}}\PD$ is $B$.

(ii) Note that $\PD$ may be an arbitrary polydisk with center at $0$;
still, Theorem \ref{ct} remains valid.
We use the unit polydisk just for simplicity.
\end{rmk}

\subsection{Proof of Behnke-Stein's Theorem}
\subsubsection{Proof of Lemma \ref{relcpt}}
{\bf (a)}
We take a subdomain $\tilde{\Omega}$ of $X$ such that
$\Omega \Subset \tilde{\Omega} \Subset X$.
Let $c \in \del \Omega$ be any point, and take
a local coordinate neighborhood system $(W_0, w)$
in $\tilde{\Omega}$ with holomorphic coordinate
$w$ such that $w=0$ at $c$.
We consider Cousin I distributions for $k=1,2, \ldots$:
\begin{align*}
\frac{1}{w^k} \quad & \hbox{on} \quad W_0, \\
0~ \quad & \hbox{on} \quad W_1=\tilde{\Omega}\setminus \{c\}.
\end{align*}
These induce cohomology classes
\[
 \left[\frac{1}{w^k}\right] \in H^1(\{W_0, W_1\}, \O_{\tilde{\Omega}})
\hookrightarrow H^1(\tilde{\Omega}, \O_{\tilde{\Omega}}),\quad
k=1,2, \ldots .
\]
Since $\dim  H^1(\tilde{\Omega}, \O_{\tilde{\Omega}})< \infty$
by \eqref{gra} (Grauert's Theorem), there is a non-trivial linear
relation over $\C$
\[
 \sum_{k=1}^\nu \gamma_k  \left[\frac{1}{w^k}\right]=0
\in H^1(\tilde{\Omega}, \O_{\tilde{\Omega}}), \quad
\gamma_k \in \C, ~ \gamma_\nu \not= 0 .
\]
Hence there is a meromorphic function $F$ on $\tilde{\Omega}$
with a pole only at $c$ such that about $c$
\begin{equation}
\label{2.1}
F(w)= \frac{\gamma_\nu}{w^\nu}+ \cdots+
\frac{\gamma_1}{w}+ \hbox{holomorphic term}.
\end{equation}
Therefore the restriction $F|_\Omega$ of $F$ to $\Omega$
is holomorphic and $\lim_{x \to c} |F(x)|=\infty$.
Thus we see that $\Omega$ is holomorphically convex.

{\bf (b)} We show the holomorphic separation property of $\Omega$
(Definition \ref{stdef} (ii)).
Let $a, b \in \Omega$ be any distinct points.
Let $F$ be the one obtained in (a) above.
If $F(a) \not= F(b)$, then it is done.
Suppose that $F(a) = F(b)$. We may assume that $F(a) = F(b)=0$.
Let $(U_0, z)$ be a local holomorphic coordinate system about
$a$ with $z(a)=0$. Then we have
\begin{equation}
\label{2.2}
F(z)=a_{k_0}z^{k_0}+ \hbox{higher order terms},
\quad a_{k_0} \not= 0, ~ k_0 \in \N,
\end{equation}
where $\N$ denotes the set of natural numbers (positive integers).
We define Cousin I distributions by
\begin{align*}
\frac{1}{z^{kk_0}} \quad & \hbox{on} \quad U_0,~k \in \N, \\
0~ \quad & \hbox{on} \quad U_1={\Omega}\setminus \{a\},
\end{align*}
which lead cohomology classes
\begin{equation}
\label{2.3}
 \left[\frac{1}{z^{kk_0}}\right] \in H^1(\{U_0, U_1\}, \O_\Omega)
\hookrightarrow H^1(\Omega, \O_\Omega),\quad
k=1,2, \ldots .
\end{equation}
It follows from \eqref{gra} that there is a non-trivial linear relation
\[
 \sum_{k=1}^\mu \alpha_k  \left[\frac{1}{z^{kk_0}}\right]=0,
\quad \alpha_k \in \C, ~ \alpha_\mu \not= 0.
\]
It follows that there is a meromorphic function $G$ on $\Omega$
with a pole only at $a$, where $G$ is written as
\begin{equation}
\label{2.4}
G(z)= \frac{\alpha_\mu}{z^{\mu k_0}}+ \cdots+
\frac{\alpha_1}{z^{k_0}}+ \hbox{holomorphic term}.
\end{equation}
With $g=G\cdot F^\mu$ we have $g \in \O(\Omega)$ and by
\eqref{2.2} and \eqref{2.4} we see that
\[
 g(a)=\alpha_\mu a_{k_0}^\mu \not= 0, \quad
g(b)=0.
\]

{\bf (c)} Let $a \in \Omega$ be any point.
We show the existence of an $h \in \O(\Omega)$ with
non vanishing differential $dh(a)\not=0$
(Definition \ref{stdef} (iii)).
Let $(U_0,z)$ be a holomorphic local coordinate system about $a$
with $z(a)=0$.
As in \eqref{2.3} we consider
\begin{equation}
\label{2.5}
 \left[\frac{1}{z^{kk_0-1}}\right] \in H^1(\{U_0, U_1\}, \O_\Omega)
\hookrightarrow H^1(\Omega, \O_\Omega),\quad
k=1,2, \ldots .
\end{equation}
In the same as above we deduce that
there is a meromorphic function $H$ on $\Omega$
with a pole only at $a$, where $H$ is written as
\begin{equation}
\label{2.6}
H(z)= \frac{\beta_\lambda}{z^{\lambda k_0-1}}+ \cdots+
\frac{\beta_1}{z^{k_0-1}}+ \hbox{holomorphic term},
\quad \beta_k \in \C, ~ \beta_{\lambda}\not=0, ~ \lambda \in \N.
\end{equation}
With $h=H\cdot F^\lambda$ we have $h \in \O(\Omega)$ and by
\eqref{2.2} and \eqref{2.6} we get
\[
d h(a)=\beta_\lambda a_{k_0}^\lambda \not= 0.
\]

Thus, $\Omega$ is Stein. \qed

\subsubsection{Proof of Lemma  \ref{rung1}}
We take a domain $\tilde{\Omega} \Subset X$ with $\tilde{\Omega} \Supset
\Omega$.
By Lemma \ref{relcpt}, $\tilde{\Omega}$ is Stein,
and hence there is a holomorphic $1$-from on $\tilde{\Omega}$
without zeros. Then we define
$ \rho(a, \Omega)$ as in \eqref{convrad} with $X=\tilde{\Omega}$.
With this $ \rho(a, \Omega)$ we have by \eqref{ct2}:
\begin{lem}
\label{hull1}
For a compact subset $K \Subset \Omega$ we get
\begin{equation*}
\rho(K, \Omega)=\rho(\hat{K}_\Omega, \Omega).
\end{equation*}
\end{lem}

\begin{lem}
\label{hull2}
Let $\Omega'$ be a domain such that $\Omega \Subset \Omega' \Subset
\tilde{\Omega}$. Assume that
\begin{equation}
\label{bcond1}
\max_{b \in \del\Omega} \rho(b, \Omega') < \rho(K, \Omega).
\end{equation}
Then,
\[
 \hat{K}_{\Omega'} \cap \Omega \Subset \Omega.
\]
\end{lem}
\begin{pf}
Since $\hat{K}_{\Omega'}$ is compact in $\Omega'$
by Lemma \ref{relcpt}, it suffices to show that
\[
 \hat{K}_{\Omega'} \cap \del \Omega = \emptyset .
\]
Suppose that there is a point
 $b \in \hat{K}_{\Omega'} \cap \del\Omega$.
It follows from Lemma \ref{hull1} that
\[
 \rho(b, \Omega') \geq \rho(\hat{K}_{\Omega'}, \Omega')
=\rho(K, \Omega') \geq \rho(K, \Omega).
\]
By assumption, $\rho(b, \Omega') < \rho(K, \Omega)$;
this is absurd.
\end{pf}

\medskip
{\bf Proof of Lemma \ref{rung1}:}
Here we use Oka's J\^oku-Ik\^o (transform to a higher space),
which is a principal method of K. Oka to reduce a difficult problem
to the one over a simpler space such as a polydisk, but of higher
dimension, and to solve it (cf.\ K. Oka \cite{ok}, e.g., \cite{nog1}).

By Lemma \ref{rung1} there are holomorphic functions
$g_j \in\O(\Omega')$ such that a finite
union $P$, called an analytic polyhedron,
 of relatively compact components
of
\[
\{x \in \Omega': |g_j(x)|<1\}
\]
satisfies ``$\hat{K}_{\Omega'} \cap \Omega \Subset P \Subset \Omega$''
and the Oka map
\[
 \Psi: x \in P \lto (g_1(x), \ldots, g_N(x)) \in \PD_N
\]
is a closed embedding into the $N$-dimensional
unit polydisk $\PD_N$.

Let $f \in \O(\Omega)$.
We identify $P$ with the image $\Psi(P) \subset \PD_N$
and regard $f|_P$ as a holomorphic function on
$\Psi(P)$.
Let $\sI$ denote the geometric ideal sheaf of the analytic
subset $\Psi(P) \subset \PD_N$.  Then we have a short
exact sequence of coherent sheaves:
\[
 0 \to \sI \to \O_{\PD_N} \to \O_{\PD_N}/\sI \to 0.
\]
By Oka's Fundamental Lemma, $H^1(\PD_N, \sI)=0$ (cf., e.g.,
\cite{nog1}, \S4.3),
which implies the surjection
\begin{equation}
\label{okafl}
 H^0(\PD_N, \O_{\PD_N}) \to H^0(\PD_N, \O_{\PD_N}/\sI)\iso \O(P) \to 0.
\end{equation}
Since $f|_P \in \O(P)$, there is an element $F \in \O(\PD_N)$
with $F|_{P}=f|_P$.
We then expand $F$ to a power series
\[
 F(w_1, \ldots, w_N)= \sum_{\nu} c_\nu w^\nu,
\quad w \in \PD_N,
\]
where $\nu$ denote multi-indices in $\{1, \ldots, N\}$.
For every $\epsilon>0$ there is a number $l \in \N$ such that
\[
\left|F(w) - \sum_{|\nu| \leq l}  c_\nu w^\nu
\right| < \epsilon,
\quad w \in \Psi(K).
\]
Substituting $w_j=g_j$, we have that
\begin{align*}
&g(x) = \sum_{|\nu| \leq l}  c_\nu g^\nu \in \O(\Omega'), \\
&|f(x) - g(x)|< \epsilon, \quad \forall x \in K.
\end{align*}\qed

\subsubsection{Proof of Theorem \ref{rung3}}\label{prung3}
We take a continuous exhaustion family $\{\Omega_t\}_{0 \leq t \leq 1}$
of subdomains of $\tilde{\Omega}$ (cf.\ Definition \ref{cex}) with
$\Omega_0=\Omega$.
Let $K \Subset \Omega$ be a compact subset and let $f \in \O(\Omega)$.
We set
\[
 T=\{t: 0 <  t \leq 1,~ \O(\Omega_t)|K 
~ \hbox{is dense in}~ \O(\Omega)|K \},
\]
where ``dense'' is taken in the sense of the maximum norm on $K$.
Note that
\begin{enumerate}
\setlength{\itemsep}{-2pt}
\item
$\rho(a, \Omega_t)$ is continuous in $t$;
\item
$\rho(K, \Omega) \leq \rho(K, \Omega_s) < \rho(K, \Omega_t)$ for
$ s<t$;
\item
$\lim_{t \searrow s } \max_{b \in \del \Omega_s} \rho(b, \Omega_t)=0$.
\end{enumerate}
It follows from Lemma \ref{rung1} that $T$ is non-empty,
open and closed. Therfore $T \ni 1$, so that
$\O(\tilde{\Omega})|K$ is dense in $\O(\Omega)|K$.
\qed

\subsubsection{Proof of Theorem \ref{bs}}
We owe the second countability axiom for Riemann surface $X$ to T. Rad\'o.
We take an increasing sequence of relatively compact domains
$\Omega_j \Subset \Omega_{j+1} \Subset X$, $j \in \N$, such that
$X=\bigcup_{j=1}^\infty \Omega_j$ and no component of
 $\Omega_{j+1} \setminus \bar\Omega_j$ is relatively
compact in $\Omega_{j+1}$.
Then, $(\Omega_j, \Omega_{j+1})$ forms a so-called Rung pair
(Theorem \ref{rung3}).
Since every $\Omega_j$ is Stein (Lemma \ref{relcpt}),
the Steiness of $X$ is deduced.
\qed

\subsection{Proofs for Riemann domains}
\subsubsection{Proof of Theorem \ref{relcomp}}
(i) Suppose that $\Omega (\Subset X)$ is a domain of holomorphy.
It follows from the assumption and Corollary \ref{dhol} that
$\Omega$ is {\it K}-complete in the sense of Grauert and
holomorphically convex. Thus, by Grauert's Theorem (\cite{g1}),
$\Omega$ is Stein.

(ii) Let $Z=\{\det d\pi=0\}$. Then, $Z$ is a thin
analytic subset of $X$.

We first take a Stein subdomain $\Omega \Subset X$
and show the plurisubharmonicity of $-\log \rho(a, \Omega)$.
 By Grauert-Remmert \cite{gr1}
it suffices to show that $-\log\rho(a, \Omega)$
is plurisubharmonic in $\Omega\setminus Z$.
Take an arbitrary point $a \in \Omega \setminus Z$,
and a complex affine line $\Lambda \subset \C^n$
passing through $\pi(a)$.
Let $\tilde\Lambda$ be the irreducible component of
$\pi^{-1}\Lambda \cap \Omega$ containing $a$.
Let $\Delta$ be a small disk about $\pi(a)$
such that $\tilde{\Delta}=\pi^{-1}\Delta \cap \tilde{\Lambda}
\Subset \tilde{\Lambda} \setminus Z$.

{\em Claim.}  The restriction
$-\log \rho(x, \Omega)|_{\tilde{\Lambda}\setminus Z}$ is subharmonic.

By a standard argument (cf.\, e.g., \cite{ho}, Proof of
Theorem 2.6.7) it suffices to prove that
if a holomorphic function $g \in \O(\tilde{\Lambda})$
satisfies
\[
 -\log \rho(x, \Omega) \leq \Re g(x), \quad x \in \del \tilde{\Delta},
\]
then
\begin{equation}
\label{inside}
 -\log \rho(x, \Omega) \leq \Re g(x), \quad x  \in \tilde{\Delta},
\end{equation}
where $\Re$ denotes the real part.
Now, we have that
\[
 \rho(x, \Omega) \geq |e^{g(x)}|, \quad x \in \del \tilde{\Delta},
\]
Since $\Omega$ is Stein, there is a holomorphic function
$f \in \O(\Omega)$ with $f|_{\tilde{\Lambda}}=g$ (cf.\ 
the arguments for \eqref{okafl}).
 Then,
\[
  \rho(x, \Omega) \geq |e^{f(x)}|, \quad x \in \del \tilde{\Delta},
\]
Since $\widehat{\tilde{\Delta}}_\Omega=\bar{\tilde{\Delta}}$,
it follows from \eqref{ct1} that
\[
  \rho(x, \Omega) \geq |e^{f(x)}|=|e^{g(x)}|, \quad x \in  \tilde{\Delta},
\]
so that \eqref{inside} follows.

Let
$\{\Omega_\nu\}_{n=1}^\infty$ be a sequence of Stein domains of $X$
such that $\Omega_\nu \Subset \Omega_{\nu+1}$ for all $\nu$
and $X=\bigcup_\nu \Omega_\nu$.
Then, $-\log \rho(a, \Omega_\nu)$, $\nu=1,2, \ldots$,
are plurisubharmonic and monotone decreasingly converges to
$-\log \rho(a, X)$. Therefore if $-\log \rho(a, X)$
is either identically $-\infty$, or plurisubharmonic ($\not\equiv - \infty$).

Suppose that $-\log \rho(a, X) \not\equiv - \infty$.
Then, the subset $A:=\{a \in X: -\log \rho(a, X) \not= -\infty$
is dense in $X$. Take any point $a \in X$ and
$U_0 (\cong \PD(\rho_0))$ as in \eqref{locbih}.
Then, there is a point $b \in A$.
Since $\rho(b, X)< \infty$, we infer that $\rho(a, X)< \infty$.
Therefore, $A=X$, and \eqref{rcont} remains valid
for $\Omega=X$. Thus, $\rho(a, X)$ is continuous in $X$.
\qed

\begin{cor}
\label{rhoplsbh}
Let $X$ be a Stein manifold satisfying Cond\,A.
Then, $-\log\rho(a, X)$ is either identically $-\infty$ or
continuous plurisubharmonic.
\end{cor}
\begin{pf}
Since $X$ is Stein, there is a holomorphic map $\pi: X \to \C^n$
which forms a Riemann domain. The assertion is immediate from
(ii) above.
\end{pf}

\begin{rmk}\rm
As a consequence, one sees with the notation
 in Corollary \ref{rhoplsbh} that
 {\em if $\Omega \subset X$ is
a domain of holomorphy, then Hartogs' radius $\rho_n(a, \Omega)$
(cf.\ Remark \ref{hradius} (ii)) is plurisubharmonic.}
This is, however, opposite to the history: The plurisubharmonicity
or the pseudoconvexity of Hartogs' radius $\rho_n(a, \Omega)$ was
found first through the study of the maximal convergence domain
of a power series (Hartogs' series)
 in several complex variables (cf.\ Oka VI \cite{ok6}, IX \cite{ok9},
Nishino \cite{nis}, Chap.\ I, Fritzsche-Grauert \cite{f-g},
 Chap.\ II).
\end{rmk}

\begin{rmk}\rm
We here give a proof of Theorem \ref{el-annar} under Cond\,A
by making use of $\rho(a,\Omega)$.

Since $\omega$ is defined in a neighborhood of $\bar{\Omega}$,
Cond\,B is satisfied at every point of the boundary $\del\Omega$;
that is, for every $b \in \del \Omega$ there are neighborhoods
$U' \Subset U \Subset X$ of $b$ such that
\[
 \rho(a, \Omega)=\rho(a, U\cap \Omega), \quad
a \in U'.
\]
If $U\cap\Omega$ is Stein, then $-\log\rho(a, \Omega)$ is 
 plurisubharmonic in $a \in U'$ by Theorem \ref{relcomp} (iii).
Therefore there is a neighborhood $V$ of $\del \Omega$ in $X$
such that $-\log\rho(a, \Omega)$ is 
 plurisubharmonic in $a \in V \cap \Omega$.
Take a real constant $C$ such that
\[
 -\log\rho(a, \Omega) < C, \quad
a \in \Omega \setminus V.
\]
Set
\[
 \psi(a)= \max\{ -\log\rho(a, \Omega), C\},
\quad a \in \Omega.
\]
Then, $\psi$ is a continuous plurisubharmonic exhaustion function
on $\Omega$. By Andreotti-Narasimhan's Theorem \ref{an-na}
in below, $\Omega$ is Stein. \qed
\end{rmk}

\subsubsection{Proof of Theorem \ref{exh1}}
In the same way as Lemma \ref{rung1} and its proof we have
\begin{lem}\label{rung2}
Let $\pi:\tilde{\Omega} \to \C^n$ be a Riemann domain such that
$\tilde{\Omega}$ satisfies Cond\,A.
Let  $\Omega \Subset \Omega' \Subset\tilde{\Omega}$ be
domains such that
\begin{equation}
\label{bcond}
\max_{b \in \del\Omega} \rho(b, \Omega') < \rho(K, \Omega).
\end{equation}
Then, every $f \in \O(\Omega)$ can be approximated uniformly
on $K$ by elements of $\O(\Omega')$.
\end{lem}

For the proof of the theorem it suffices to show that
$(\Omega_t, \Omega_s)$ is a Runge pair for $0 \leq t < s <1$.
Since any fixed $\Omega_{s'}$ ($s<s'<1$) satisfies Cond\,A, we have
the scalar $\rho(a, \Omega_s)$.
Take a compact subset $K \Subset \Omega_t$.
Then, for $s>t$ sufficiently close to $t$ we have
\[
\max_{b \in \del\Omega_t} \rho(b, \Omega_s) < \rho(K, \Omega_t).
\]
It follows from Lemma \ref{rung2} that 
$\O(\Omega_s)|_K$ is dense if $\O(\Omega_t)|_K$. 
Then, the rest of the proof is the same as in \S\ref{prung3}.
\qed

\subsubsection{Proof of Theorem \ref{levi2}}
Here we will use the following result:
\begin{thm}[Andreotti-Narasimhan \cite{an}]
\label{an-na}
Let $\pi: X \to \C^n$ be a Riemann domain.
If $X$ admits a continuous plurisubharmonic exhaustion function,
then $X$ is Stein.
\end{thm}

For an ideal boundary point $b \in \del X$ there are
connected open subsets $\tilde{V} \subset \widetilde{W}$ as
in Cond\,B such that
\begin{equation}
\label{lp1}
\rho(a, X) = \rho(a, \widetilde{W}).
\end{equation}
By the assumption, $\widetilde{W}$ can be chosen to be Stein.
By Theorem \ref{relcomp} (ii), $-\log \rho(a, \widetilde{W})$
 is plurisubharmonic in $a \in \widetilde{V}$, and hence
so is $-\log \rho(a, X)$ in $\widetilde{V}$.
Since $\lim_{a \to \del X}\rho(a, X)=0$ by Cond\,B,
 there is a closed subset $F \subset X$ such that
\begin{enumerate}
\item
$F \cap \{x \in X: \|\pi(x)\|\leq R\}$
is compact for every $R>0$,
\item
$ -\log\rho(a, X)$ is plurisubharmonic in $a \in X\setminus F$,
\item
$ -\log\rho(a, X) \to \infty$ as $a \to \del X$.
\end{enumerate}
From this we may construct a continuous plurisubharmonic
exhaustion function on $X$ as follows:

We fix a point $a_0 \in F$, and may assume that
$\pi(a_0)=0$. Let $X_\nu$ be a connected component
of $\{ \|\pi\|<\nu \}$ containing $a_0$. Then,
$\bigcup_\nu X_\nu=X$. Put
\[
 \Omega_\nu=X_\nu \setminus F \Subset X.
\]
Take a real constant $C_1$ such that
\[
-\log \rho(a, X) < C_1, \quad a \in \bar{\Omega}_1.
\]
Then we set
\[
 \psi_1(a)=\max\{-\log \rho(a, X),  C_1\},
\quad a \in X.
\]
Then, $\psi_1$ is plurisubharmonic in $X_1$.
We take a positive constant $C_2$ such that
\[
 -\log \rho(a, X) < C_1+C_2 (\|\pi(a)\|^2-1)^+, \quad
 a \in \bar{\Omega}_2,
\]
where $(\cdot)^+=\max\{\cdot, 0\}$.
Put
\begin{align*}
p_2(a) &= C_1+C_2 (\|\pi(a)\|^2-1)^+, \\
\psi_2(a) &=\max\{-\log \rho(a, X), p_2(a)\}.
\end{align*}
Then $\psi_1(a)=\psi_2(a)$ in $a \in X_1$ and $\psi_2(a)$
is plurisubharmonic in $X_2$.
Similarly, we take $C_3 > C_2$ so that
\[
 -\log \rho(a, X) < p_2(a) +C_3 (\|\pi(a)\|^2-2^2)^+, \quad
 a \in \bar{\Omega}_3,
\]
Put
\begin{align*}
p_3(a) &= p_2(a)+C_3 (\|\pi(a)\|^2-2^2)^+, \\
\psi_3(a) &=\max\{-\log \rho(a, X), p_3(a)\}.
\end{align*}
Then $\psi_3(a)=\psi_2(a)$ in $a \in X_2$ and $\psi_3(a)$
is plurisubharmonic in $X_3$.
Inductively, we may take a continuous function
$\psi_\nu(a)$, $\nu=1,2, \ldots$,
 such that $\psi_\nu$ is plurisubharmonic
in $X_\nu$ and $\psi_{\nu+1}|_{X_\nu}=\psi_\nu|_{X_\nu}$.
it is clear from the construction that
\[
 \psi(a)=\lim_{\nu \to \infty} \psi_\nu(a), \quad a \in X
\]
is a continuous plurisubharmonic exhaustion function of $X$.

Finally by Andreotti-Narasimhan's Theorem \ref{an-na} we see that
$X$ is Stein.

\section{Examples and some more on $\rho(a, X)$}
\label{smr}

{\bf (a)} (Grauert's example).
Grauert \cite{nar} gave a counter-example to the Levi problem
for ramified Riemann domains over $\pnc$:
There is a locally Stein domain $\Omega$ in a complex torus
$M$ such that $\O(\Omega)=\C$.
Then, $M$ satisfies Cond\,A.
 One may assume that $M$ is projective
algebraic, so that there is a holomorphic finite map
$\tilde{\pi}: M \to \pnc$, which is a Riemann domain over $\pnc$.
Then, the restriction $\pi=\tilde{\pi}|_\Omega: \Omega \to \pnc$
is a Riemann domain over $\pnc$, which satisfies Cond\,A and Cond\,B.
Therefore, Theorem \ref{levi2} cannot be extended to
a Riemann domain over $\pnc$.

\begin{rmk}\rm
Let $\pi: \Omega \to \pnc$ be Grauert's example as above.
Let $\C^n$ be an affine open subset of $\pnc$, and
let $\pi': \Omega' \to \C^n$ be the restriction of
$\pi: \Omega \to \pnc$ to $\C^n$.
Then, {\em  $\Omega'$ is Stein by Theorem \ref{levi2}.}

The Steinness of $\Omega'$ may be not inferred
by a formal combination of the known results on
pseudoconvexity, since it is an unbounded domain
(cf., e.g., \cite{nar}, \cite{lieb}).
\end{rmk}

\smallskip
{\bf (b)}
Domains in the products of open Riemann surfaces and complex semi-tori
(cf.\ \cite{nw}, Chap.\ 5)
serve for examples satisfying Cond\,A.

\smallskip
{\bf (c)}
An open Riemann surface $X$ is {\em not} Kobayashi hyperbolic
if and only if $X$ is biholomorphic to $\C$ or $\C^*=\C\setminus\{0\}$
(For the Kobayashi hyperbolicity in general, cf.\ \cite{ko}, \cite{nw}).

{\bf (c1)} Let $X=\C$. If $\omega=dz$, then $\rho(a, dz, \C) \equiv \infty$
for every $a \in \C$. If $\omega=e^zdz$, then a simple calculation
implies that
\[
 \rho(a, e^zdz, \C)=|e^a|.
\]

{\bf (c2)} Let $X=\C^*$. If $\omega=z^k dz$ with $k \in \Z\setminus\{-1\}$,
 then
$$ \rho(a, z^k dz, \C^*)=\left| \frac{1}{k+1} a^{k+1} \right|.
$$
Therefore, $\lim_{a \to 0} \rho(a, z^k dz, \C^*)=0$ for $k \geq 0$,
and $\lim_{a \to \infty} \rho(a, z^k dz, \C^*)=0$ for $k \leq -2$.
If $\omega=\frac{dz}{z}$, then $\rho(a,\frac{dz}{z}, \C^*) \equiv
 \infty$.
It follows that
$$
\psi(a) := \max\{-\log \rho(a, dz, \C^*),
-\log  \rho(a, z^{-2} dz, \C^*) \}
$$
is continuous subharmonic in $\C^*$, and
$\lim_{a \to 0, \infty} \psi(a)=\infty$.

Thus, the finiteness or the infiniteness of $ \rho(a, \omega, X)$
depends on the choice of $\omega$.

{\bf (d)} For a Kobayashi hyperbolic open Riemann surface $X$ we take a
holomorphic $1$-form $\omega$ without zeros, and write
\[
    \|\omega(a)\|_{X}=|\omega(v)|, \quad v \in \mathbf{T}(X)_a,
~ F_X(v)=1.
   \]
Then it follows from \eqref{kh} that
 $\rho(a, \omega, X)\leq \|\omega(a)\|_{X}$.
We set
\begin{align*}
\rho^+(a,X) &= \sup \{\rho(a, \omega, X): \omega ~
\hbox{hol.\ 1-form without zeros},~\|\omega(a)\|_{X}=1\},\\
\rho^-(a, X) &= \inf \{\rho(a, \omega, X): \omega ~
\hbox{hol.\ 1-form without zeros}
,~\|\omega(a)\|_{X}=1\}.
\end{align*}
Clearly, $\rho^{\pm}(a, X) (\leq 1)$ are biholomorphic invariants
of $X$, but we do not know the behavior of them.

\rightline{Graduate School of Mathematical Sciences}
\rightline{University of Tokyo}
\rightline{Komaba, Meguro-ku, Tokyo 153-8914}
\rightline{Japan}
\rightline{e-mail: noguchi@ms.u-tokyo.ac.jp}

\end{document}